\let\ReportLabel\label
\let\ReportRefstepcounter\refstepcounter
\AddToHook{package/hyperref/before}{%
  \let\label\ReportLabel
  \let\refstepcounter\ReportRefstepcounter}
\documentclass{siamart220329}
\usepackage{amsmath,amssymb,bm}
\usepackage{graphicx}
\usepackage{algorithmic}
\usepackage{placeins}
\usepackage{enumitem}
\usepackage{tabularx,booktabs}
\usepackage{xurl}
\usepackage[numbers]{natbib}
\usepackage{etoolbox}

\setlist[enumerate]{leftmargin=2em}
\setlist[itemize]{leftmargin=2em}

\newsiamremark{remark}{Remark}
\newsiamthm{assumption}{Assumption}
\DeclareMathOperator*{\argmin}{arg\,min}

\hypersetup{pdftitle={Recent Developments in HPR Methods for Large-Scale Linear Programming}}

\newcommand{\ReportBlackBibitem}[2][]{%
  \par\color{black}%
  \hypersetup{linkcolor=\ReportSavedLinkColor,
    citecolor=\ReportSavedCiteColor,urlcolor=\ReportSavedUrlColor}%
  \ifboolexpr{test {\ifstrequal{#2}{lu2023cupdlp}} or
    test {\ifstrequal{#2}{optbin-pdlp-large}}}
    {\color{black}\hypersetup{linkcolor=black,citecolor=black,urlcolor=black}}{}%
  \ReportSavedBibitem[#1]{#2}}
\makeatletter
\AtBeginEnvironment{thebibliography}{%
  \let\ReportSavedBibitem\bibitem
  \let\ReportSavedLinkColor\@linkcolor
  \let\ReportSavedCiteColor\@citecolor
  \let\ReportSavedUrlColor\@urlcolor
  \let\bibitem\ReportBlackBibitem}
\makeatother

\title{Recent Developments in HPR Methods for  \\ Large-Scale Linear Programming}
\date{1 November 2026}

\makeatletter
\patchcmd{\@maketitle}
  {\bfseries\MakeUppercase{\siampretitle}\@ucnothanks{\@title}}
  {\large\bfseries\siampretitle\@title}
  {}{\PackageError{hpr-lp}{Unable to format the title}{Check the SIAM title definition.}}
\patchcmd{\@maketitle}
  {\footnotesize
      \spaceskip=1.5\fontdimen2\the\and{\ \scriptsize AND }
      \@tempswafalse\expandafter\@doauth\@author\and\relax}
  {\small\@author\par}
  {}{\PackageError{hpr-lp}{Unable to format the author block}{Check the SIAM title definition.}}
\patchcmd{\@maketitle}
  {\@subtitle}
  {\@subtitle\par\vspace{0.6em}{\small\@date\par}}
  {}{\PackageError{hpr-lp}{Unable to format the date}{Check the SIAM title definition.}}
\makeatother
\newcommand{\authoremail}[1]{\href{mailto:#1}{\texttt{#1}}}
\author{%
Kaihuang Chen\textsuperscript{1}\quad
Qilin Lai\textsuperscript{1}\quad
Yijia Liang\textsuperscript{1}\quad
Jiawei Lu\textsuperscript{1}\quad
Defeng Sun\textsuperscript{1}\\
Yancheng Yuan\textsuperscript{1}\quad
Guojun Zhang\textsuperscript{1}\quad
Xiangyu Zhang\textsuperscript{1}\quad
Xinyuan Zhao\textsuperscript{2}}
\subtitle{%
\textsuperscript{1}Department of Applied Mathematics, The Hong Kong Polytechnic University,\par
Hung Hom, Hong Kong\par
\textsuperscript{2}Department of Mathematics, Beijing University of Technology, Beijing, P.R. China\par\vspace{0.3em}
\authoremail{kaihuang.chen@connect.polyu.hk}\quad \authoremail{jiawei.lu@connect.polyu.hk}\par
\authoremail{defeng.sun@polyu.edu.hk}\quad \authoremail{yancheng.yuan@polyu.edu.hk}\par
\authoremail{guojun.zhang@connect.polyu.hk}\quad \authoremail{xiangyu.zhang@connect.polyu.hk}\par
\authoremail{xyzhao@bjut.edu.cn}
}
\hypersetup{pdfauthor={Kaihuang Chen, Qilin Lai, Yijia Liang, Jiawei Lu, Defeng Sun, Yancheng Yuan, Guojun Zhang, Xiangyu Zhang, Xinyuan Zhao}}

\begin{document}
\maketitle
\begin{abstract}
We present recent developments in Halpern Peaceman--Rachford (HPR)
methods for solving large-scale linear programming problems.
We first introduce the semi-proximal HPR framework and summarize
its convergence guarantees, together with its connections to
Halpern PDHG and ergodic Peaceman--Rachford iterations.
To improve practical performance, we develop adaptive restart
and penalty parameter update strategies.
We then describe several computational enhancements,
including GPU-based preprocessing, diagonal scaling,
reusable execution plans for sparse matrix--vector multiplication,
and sieving techniques.
We further discuss multi-GPU parallelization and lossless
model storage and loading for handling large-scale instances.
These developments form the foundation of HPR-LP,
a GPU-accelerated solver for large-scale linear programming.
\end{abstract}
\begin{keywords}
linear programming, Halpern Peaceman--Rachford method, GPU computing,
preprocessing, scaling, sieving.
\end{keywords}
\begin{MSCcodes}
90C05, 90C25, 65Y20
\end{MSCcodes}

\bigskip
\noindent{\large\bfseries Outline}\par\nobreak\medskip
\begingroup
\newcommand{\outlineentry}[1]{%
  \hyperref[#1]{\ref*{#1}} & \hyperref[#1]{\nameref*{#1}} & \pageref*{#1}\\[0.5em]}
\noindent\begin{tabularx}{\textwidth}{@{}l@{\hspace{0.8em}}>{\raggedright\arraybackslash}Xr@{}}
\outlineentry{sec:intro}
\outlineentry{sec:2}
\outlineentry{sec:3}
\outlineentry{sec:implementation}
\outlineentry{sec:preprocessing}
\outlineentry{sec:scaling}
\outlineentry{sec:spmv}
\outlineentry{sec:sieving}
\outlineentry{sec:multigpu}
\outlineentry{sec:optbin}
\outlineentry{sec:conclusion}
\end{tabularx}
\endgroup
\clearpage

\section{Introduction}\label{sec:intro}

\href{https://github.com/PolyU-IOR/HPR-LP}{HPR-LP} applies the Halpern Peaceman--Rachford (HPR) method to large-scale
linear programming (LP). Its basic iteration combines semi-proximal
splitting with Halpern anchoring and can be implemented using sparse
matrix--vector products and box projections. The method and its original
implementation are developed in~\cite{sun2025accelerating,chen2024hpr}.
This report brings together the algorithm and subsequent implementation
techniques in a common notation. We emphasize the role of each component
and refer to the original sources for established analyses and experiments.

\subsection{Problem formulation}
We consider
\begin{equation}\label{model:primalLP}
 \min_{x\in\mathbb R^n}\ \langle c,x\rangle
 \quad\text{subject to}\quad Ax\in\mathcal K,\qquad x\in\mathcal C,
\end{equation}
where $A\in\mathbb R^{m\times n}$, $c\in\mathbb R^n$, and
\[
 \mathcal K=\{s\in\mathbb R^m:l_c\leq s\leq u_c\},\qquad
 \mathcal C=\{x\in\mathbb R^n:l_v\leq x\leq u_v\}.
\]
Lower bounds may be $-\infty$ and upper bounds may be $+\infty$,
with $l_c\leq u_c$ and $l_v\leq u_v$.
Equal lower and upper bounds represent equalities or fixed variables.
Writing the dual in minimization form gives
\begin{equation}\label{model:dualLP}
 \min_{y\in\mathbb R^m,\,z\in\mathbb R^n}
 \ \delta_{\mathcal K}^*(-y)+\delta_{\mathcal C}^*(-z)
 \quad\text{subject to}\quad A^*y+z=c.
\end{equation}
Here $\delta_S$ is the indicator of a closed convex set $S$ and
$\delta_S^*$ is its support function. The dual objective in this convention
is the negative of the usual maximization objective; primal and dual
optimal objective values therefore sum to zero.

\subsection{Notation and scope}
We use $A^*$ for the transpose and $\langle\cdot,\cdot\rangle$ for the
Euclidean inner product. The norm $\|\cdot\|$ is Euclidean for vectors
and spectral for linear operators. For a self-adjoint positive
semidefinite operator $\mathcal M$,
$\|v\|_{\mathcal M}^2=\langle v,\mathcal Mv\rangle$.
The symbols $\partial f$, $\Pi_S$, and $\mathcal N_S$ denote the convex
subdifferential, Euclidean projection, and normal cone, respectively.
We write $w=(y,z,x)\in\mathbb W:=\mathbb R^m\times\mathbb R^n\times\mathbb R^n$.
For iteration variables, a bar denotes the output of one semi-proximal
sweep, a hat its reflected point, and a subscript $a$ an ergodic average.
The indices $k$, $r$, and
$t$ denote an unrestarted iteration, a restart cycle, and an iteration
within a cycle. In restarted algorithms, $k$ also counts total iterations.
The abbreviation SpMV means sparse matrix--vector multiplication.

Iteration formulas reuse this notation for the working LP after any
preprocessing and scaling; original-model checks require recovery to
the input LP.

Section~\ref{sec:2} presents the HPR method for LP, and
Section~\ref{sec:3} discusses its relationships with other splitting methods.
Section~\ref{sec:implementation} describes restart and penalty updates.
Sections~\ref{sec:preprocessing} and~\ref{sec:scaling} address preprocessing
and scaling before iteration.
Section~\ref{sec:spmv} describes ISpMV and its reusable execution plans;
Section~\ref{sec:sieving} presents a sieving technique;
and Section~\ref{sec:multigpu} describes multi-GPU parallelization.
Section~\ref{sec:optbin} summarizes OptBin as a separate storage and loading
component.

\subsection{Software repositories}\label{subsec:repositories}
The HPR project website is \url{https://hpr.top/}.
Table~\ref{tab:repositories} lists the public GitHub repositories for
our Julia, C/CUDA, Python, and MATLAB implementations of HPR-LP,
together with the GPU preprocessing framework.

\begin{table}[htbp]
\centering\small
\begin{tabularx}{\textwidth}{@{}ll>{\raggedright\arraybackslash}X@{}}
\toprule
Software & Execution & GitHub repository\\
\midrule
HPR-LP (Julia) & CPU, GPU &
\href{https://github.com/PolyU-IOR/HPR-LP}{\nolinkurl{PolyU-IOR/HPR-LP}}\\
HPR-LP-C (C/CUDA) & GPU &
\href{https://github.com/PolyU-IOR/HPR-LP-C}{\nolinkurl{PolyU-IOR/HPR-LP-C}}\\
HPR-LP-Python & GPU &
\href{https://github.com/PolyU-IOR/HPR-LP-Python}{\nolinkurl{PolyU-IOR/HPR-LP-Python}}\\
HPR-LP-MATLAB & CPU &
\href{https://github.com/PolyU-IOR/HPR-LP-MATLAB}{\nolinkurl{PolyU-IOR/HPR-LP-MATLAB}}\\
GPU-Preprocessor & GPU &
\href{https://github.com/PolyU-IOR/GPU-Preprocessor}{\nolinkurl{PolyU-IOR/GPU-Preprocessor}}\\
\bottomrule
\end{tabularx}
\caption{HPR-LP implementations and GPU preprocessing software.}
\label{tab:repositories}
\end{table}
\FloatBarrier

HPR-LP-C also provides a C/C++ API and
\href{https://github.com/PolyU-IOR/HPR-LP-C/blob/main/bindings/python/README.md}{Python},
\href{https://github.com/PolyU-IOR/HPR-LP-C/blob/main/bindings/julia/README.md}{Julia (including JuMP)}, and
\href{https://github.com/PolyU-IOR/HPR-LP-C/blob/main/bindings/matlab/README.md}{MATLAB}
interfaces. These interfaces call the same C/CUDA GPU solver;
the standalone Julia, Python, and MATLAB implementations have their own
repositories in Table~\ref{tab:repositories}.

\section{The HPR method for LP}\label{sec:2}

\subsection{Semi-proximal iteration}
For $\sigma>0$, define the augmented Lagrangian of~\eqref{model:dualLP} by
\begin{equation}\label{def:augmented-lagrangian}
 L_\sigma(y,z;x)=\delta_{\mathcal K}^*(-y)+\delta_{\mathcal C}^*(-z)
 +\langle x,A^*y+z-c\rangle+\frac{\sigma}{2}\|A^*y+z-c\|^2.
\end{equation}
Let $\mathcal T_1\succeq0$ be self-adjoint with
$\mathcal T_1+AA^*\succ0$. Algorithm~\ref{alg:sp-HPR} is the
semi-proximal HPR method of~\cite{sun2025accelerating,chen2024hpr}.
Its first three steps define a sweep $\mathcal D_\sigma(w^k)=\bar w^{k+1}$;
the reflected map is $\mathcal F_\sigma=2\mathcal D_\sigma-I$.
The final step applies Halpern anchoring~\cite{halpern1967fixed,lieder2021convergence}.

\begin{algorithm}[htbp]
\caption{Semi-proximal HPR for LP~\cite{sun2025accelerating,chen2024hpr}}
\label{alg:sp-HPR}
\begin{algorithmic}[1]
\REQUIRE $w^0\in\mathbb W$, $\sigma>0$, $\mathcal T_1\succeq0$,
         and $\mathcal T_1+AA^*\succ0$.
\FOR{$k=0,1,\ldots$}
 \STATE $\displaystyle\bar z^{k+1}=\argmin_z L_\sigma(y^k,z;x^k)$.
 \STATE $\displaystyle\bar x^{k+1}=x^k+\sigma(A^*y^k+\bar z^{k+1}-c)$.
 \STATE $\displaystyle\bar y^{k+1}=\argmin_y
   \{L_\sigma(y,\bar z^{k+1};\bar x^{k+1})
     +\tfrac{\sigma}{2}\|y-y^k\|_{\mathcal T_1}^2\}$.
 \STATE $\hat w^{k+1}=2\bar w^{k+1}-w^k$.
 \STATE $\displaystyle w^{k+1}=\frac{1}{k+2}w^0+\frac{k+1}{k+2}\hat w^{k+1}$.
\ENDFOR
\ENSURE The sequence $\{\bar w^{k+1}\}_{k\geq0}$.
\end{algorithmic}
\end{algorithm}

\subsection{Projection-based implementation}\label{subsec:projection}
The standard matrix-free choice is
\begin{equation}\label{eq:proximal-choice}
 \mathcal T_1=\lambda_A I_m-AA^*,\qquad
 \lambda_A>0,\quad\lambda_A\geq\|A\|^2.
\end{equation}
With this choice, the sweep has the explicit form
\begin{equation}\label{HPR-LP-2}
\begin{aligned}
 \xi^k&=x^k+\sigma(A^*y^k-c),
 &\bar x^{k+1}&=\Pi_{\mathcal C}(\xi^k),\\
 \bar z^{k+1}&=(\bar x^{k+1}-\xi^k)/\sigma,\\
 \zeta^k&=A(2\bar x^{k+1}-x^k)-\sigma\lambda_Ay^k,
 &\bar y^{k+1}&=\frac{\Pi_{\mathcal K}(\zeta^k)-\zeta^k}{\sigma\lambda_A}.
\end{aligned}
\end{equation}
Thus each sweep requires products with $A$ and $A^*$ and componentwise
operations, including two box projections. This structure is well suited
to GPU execution~\cite[Section~3.3]{chen2024hpr}.

For equality constraints $\mathcal K=\{b\}$, the alternative
$\mathcal T_1=0$ gives
\begin{equation}\label{eq:AATy=R}
 AA^*\bar y^{k+1}
 =\frac{b-A\bar x^{k+1}}{\sigma}-A(\bar z^{k+1}-c).
\end{equation}
Under the assumptions of Algorithm~\ref{alg:sp-HPR}, this choice requires
$AA^*\succ0$.

\subsection{Established convergence guarantees}
The KKT conditions for~\eqref{model:primalLP}--\eqref{model:dualLP} are
\begin{equation}\label{eq:KKT}
 -y^*\in\mathcal N_{\mathcal K}(Ax^*),\qquad
 -z^*\in\mathcal N_{\mathcal C}(x^*),\qquad A^*y^*+z^*=c.
\end{equation}
Throughout the convergence statements, we impose the standard assumption
below.
\begin{assumption}\label{ass:CQ}
The KKT system~\eqref{eq:KKT} has a solution.
\end{assumption}
Define
\begin{equation}\label{def:M}
 \mathcal M_\sigma=
 \begin{bmatrix}
 \sigma(AA^*+\mathcal T_1)&0&A\\
 0&0&0\\ A^*&0&\sigma^{-1}I_n
 \end{bmatrix},
\end{equation}
and the KKT residual
\begin{equation}\label{def:KKT_residual}
 \mathcal R(w)=\begin{pmatrix}
 Ax-\Pi_{\mathcal K}(Ax-y)\\
 x-\Pi_{\mathcal C}(x-z)\\
 c-A^*y-z
 \end{pmatrix}.
\end{equation}
The operator $\mathcal M_\sigma$ is positive semidefinite and generally
singular; $\mathcal R(w)=0$ is equivalent to~\eqref{eq:KKT}.

\begin{theorem}[\cite{sun2025accelerating}]\label{Th:complexity-acc-pADMM}
Under Assumption~\ref{ass:CQ}, with fixed parameters as in
Algorithm~\ref{alg:sp-HPR}, the sequence $\{\bar w^{k+1}\}$ converges to
a KKT point $w^*$. Set $R_0=\|w^0-w^*\|_{\mathcal M_\sigma}$ and
$d(y,z)=\delta_{\mathcal K}^*(-y)+\delta_{\mathcal C}^*(-z)$.
For every $k\geq0$,
\begin{subequations}\label{eq:complexity-bounds}
\begin{align}
 \|\bar w^{k+1}-w^k\|_{\mathcal M_\sigma}
 &\leq\frac{R_0}{k+1},\label{eq:complexity-bound-M}\\
 \|\mathcal R(\bar w^{k+1})\|
 &\leq\frac{\sigma(\|A\|+\|\sqrt{\mathcal T_1}\|)+1}{\sqrt\sigma}
       \frac{R_0}{k+1},\label{eq:complexity-bound-KKT}\\
 -\frac{\|x^*\|}{\sqrt\sigma}\frac{R_0}{k+1}
 &\leq d(\bar y^{k+1},\bar z^{k+1})-d(y^*,z^*)
 \leq\left(3R_0+\frac{\|x^*\|}{\sqrt\sigma}\right)\frac{R_0}{k+1}.
 \label{eq:complexity-bound-obj}
\end{align}
\end{subequations}
\end{theorem}
These are the LP specialization of the convergence and nonergodic
complexity results in~\cite[Corollary~3.5, Proposition~2.9, and
Theorem~3.7]{sun2025accelerating}.
The bounds concern the base method with a fixed metric and anchor.
They motivate the practical rules in Section~\ref{sec:implementation}
but do not, by themselves, establish the same global bounds for arbitrary
restarts, penalty changes, or sieved iterations.

\section{Relationships with other splitting methods}\label{sec:3}

This section summarizes the relationships established
in~\cite[Section~3]{chen2027relationships}, using the notation of this report.
The relationships concern update maps under matched parameters;
solver implementations may differ in their restarts, scaling,
stopping tests, and GPU kernels.

\subsection{Halpern PDHG and cuPDLPx}\label{subsec:pdhg}
Let $u=(y,x)$ and let $T_{\rm PDHG}$ be the PDHG sweep with primal
stepsize $\eta/\omega$ and dual stepsize $\eta\omega$:
\begin{equation}\label{cuPDLPx}
\begin{aligned}
 \bar x&=\Pi_{\mathcal C}\bigl(x+\tfrac{\eta}{\omega}(A^*y-c)\bigr),\\
 q&=A(2\bar x-x)-\tfrac{1}{\eta\omega}y,\qquad
 \bar y=\eta\omega\bigl(\Pi_{\mathcal K}(q)-q\bigr),\\
 T_{\rm PDHG}(u)&=(\bar y,\bar x).
\end{aligned}
\end{equation}
The reflected Halpern PDHG framework used by cuPDLPx applies
\begin{equation}\label{eq:relaxed-halpern}
 u^{k+1}=\frac{1}{k+2}u^0+
 \frac{k+1}{k+2}\bigl((1+\gamma)T_{\rm PDHG}(u^k)-\gamma u^k\bigr),
 \qquad 0\leq\gamma\leq1,
\end{equation}
with fixed parameters within the base iteration~\cite{lu2024restarted,lu2025cupdlpx}.
The parameter correspondence is
\begin{equation}\label{eq:sigma&lambda_A}
 \sigma=\frac{\eta}{\omega},\qquad
 \lambda_A=\frac{1}{\eta^2}\geq\|A\|^2,\qquad \eta,\omega>0.
\end{equation}
Under~\eqref{eq:sigma&lambda_A}, the sweep~\eqref{cuPDLPx} is precisely
the $(y,x)$ part of~\eqref{HPR-LP-2}. With $\gamma=1$ and the same
initial $(y^0,x^0)$, the HPR and reflected Halpern PDHG sequences
therefore agree~\cite[Proposition~3.1]{chen2027relationships}.
More generally, $\rho=1+\gamma$ gives the corresponding
Halpern-accelerated preconditioned ADMM relaxation
in~\cite{sun2025accelerating}; this connection is also discussed
in~\cite[Section~1]{chen2024hpr}.

For $A\ne0$, the endpoint $\eta=1/\|A\|$ is allowed by the
semidefinite metric analysis, since~\eqref{eq:proximal-choice} still
holds. The relaxation $\rho=2$ gives the favorable factor
$1/(\rho(k+1))$ in the cited residual bound.

\subsection{Ergodic PR and affine Halpern iterations}\label{subsec:epr}
Using the sweep $\mathcal D_\sigma$ from Section~\ref{sec:2}, the
unanchored PR iteration is
\begin{equation}\label{eq:epr}
 v^{j+1}=\mathcal F_\sigma(v^j),\qquad
 \bar v^{j+1}=\mathcal D_\sigma(v^j),\qquad v^0=w^0.
\end{equation}
We use $v$ to distinguish this trajectory from the HPR trajectory $w$.
EPR averages the sweep outputs. In the indexing used here, set
\begin{equation}\label{eq:epr-averages}
 \bar v_a^{k+1}=\frac{1}{k+2}\sum_{j=0}^{k+1}\bar v^{j+1},
 \qquad
 v_a^{k+1}=\frac{1}{k+2}\sum_{j=0}^{k+1}v^j.
\end{equation}
Convergence of the sweep-output averages $\bar v_a^{k+1}$ is established
in~\cite{chen2025peaceman}. That paper also gives an example in which
the raw PR trajectory is unbounded and its averages $v_a^{k+1}$ fail
to converge. These two averages must therefore be distinguished.

For an affine map $\mathcal F(w)=Rw+b$, Halpern iteration with anchor
$w^0$ and weight $1/(k+2)$ satisfies
\begin{equation}\label{eq:ExH-aff-avg}
 w^{k+1}=\frac{1}{k+2}\sum_{j=0}^{k+1}\mathcal F^j(w^0).
\end{equation}
This affine Halpern--ergodic identity is stated and proved
in~\cite[Proposition~3.2]{chen2027relationships}.
The linear case is discussed in~\cite[Section~1]{kornlein2015quantitative},
and the affine PDHG specialization appears
in~\cite[Proposition~1]{lu2024restarted}.

Under~\eqref{eq:proximal-choice}, the sweep uses piecewise-affine box
projections.
Fixing their lower-active, upper-active, and interior components defines
a convex region on which $\mathcal F_\sigma$ is affine. If the entire
PR trajectory~\eqref{eq:epr} stays in one such region, its convex
averages stay there as well. With the same anchor and parameters,
\cite[Corollary~3.1]{chen2027relationships} gives
\begin{equation}\label{cor:hpr-epr-equivalence}
 w^{k+1}=v_a^{k+1}.
\end{equation}
This identity concerns averages of raw PR iterates, not the
sweep-output averages whose convergence is proved in the EPR analysis.

\subsection{Identification of projection patterns}\label{subsec:identification}
The finite-identification result below uses strict complementarity.
For each nonfixed variable,
strict complementarity means $z_j^*>0$ at its lower bound,
$z_j^*<0$ at its upper bound, and $z_j^*=0$ in its interior.
For each nonequality constraint, the analogous signs are
$y_i^*>0$, $y_i^*<0$, and $y_i^*=0$ at the lower bound, upper bound,
and interior of $\mathcal K$, respectively. Equality-bound components
need no sign condition.

\begin{proposition}[{\cite[Proposition~3.3]{chen2027relationships}}]\label{thm:act-stab-bar}
Suppose Assumption~\ref{ass:CQ} and~\eqref{eq:proximal-choice} hold,
and the limit $w^*$ of the sweep outputs in Algorithm~\ref{alg:sp-HPR}
satisfies the strict complementarity conditions above.
Then the projection patterns of both $\Pi_{\mathcal C}(\xi^k)$ and
$\Pi_{\mathcal K}(\zeta^k)$ agree with those at $w^*$ for all
sufficiently large $k$.
\end{proposition}

As clarified in~\cite[Remark~3.3]{chen2027relationships}, finite
identification concerns the HPR trajectory with fixed parameters.
For~\eqref{cor:hpr-epr-equivalence}, HPR and PR must start from a common
point and the subsequent PR iterates must remain in the identified region;
finite identification alone does not imply eventual equality of two
previously generated trajectories.

\section{Restart and penalty updates}\label{sec:implementation}

HPR-LP combines the base iteration with adaptive restarts and penalty
updates~\cite[Sections~3.1--3.2]{chen2024hpr}; a unified description is
given in~\cite[Section~2.2]{chen2027relationships}. Within cycle $r$, the anchor
$w^{r,0}$ and penalty $\sigma_r$ remain fixed. We write
$\mathcal M_r=\mathcal M_{\sigma_r}$ and reuse the sweep
$\mathcal D_{\sigma_r}$ defined in Section~\ref{sec:2}.

\begin{algorithm}[htbp]
\caption{Restarted HPR-LP~\cite{chen2024hpr}}
\label{alg:HPR_LP}
\begin{algorithmic}[1]
\REQUIRE $w^{0,0}\in\mathbb W$, $\sigma_0>0$, and $\mathcal T_1$
as in Algorithm~\ref{alg:sp-HPR}.
\STATE Set $r=0$ and the total iteration counter $k=0$.
\LOOP
 \STATE Set $t=0$.
 \REPEAT
  \STATE $\bar w^{r,t+1}=\mathcal D_{\sigma_r}(w^{r,t})$.
  \STATE $\hat w^{r,t+1}=2\bar w^{r,t+1}-w^{r,t}$.
  \STATE $\displaystyle w^{r,t+1}=\frac{1}{t+2}w^{r,0}
                +\frac{t+1}{t+2}\hat w^{r,t+1}$.
  \STATE Set $t=t+1$ and $k=k+1$.
  \STATE \textbf{return} $\bar w^{r,t}$ if the termination test is satisfied.
 \UNTIL{a restart criterion in~\eqref{eq:restart-rules} holds}
 \STATE Set $\tau_r=t$ and $w^{r+1,0}=\bar w^{r,\tau_r}$.
 \STATE Update $\sigma_{r+1}$ as described in Section~\ref{subsec:sigma};
        set $r=r+1$.
\ENDLOOP
\end{algorithmic}
\end{algorithm}

\subsection{Restart criteria}\label{subsec:restart}
The distance $R_0$ in Theorem~\ref{Th:complexity-acc-pADMM} depends on
the anchor. Resetting the anchor near a solution can reduce this distance,
but it is not directly observable. HPR-LP instead monitors the sweep
residual
\begin{equation}\label{eq:restart-residual}
 \varepsilon_{r,t}
   =\|w^{r,t-1}-\bar w^{r,t}\|_{\mathcal M_r},
 \qquad t\geq1.
\end{equation}
This is one half of the corresponding reflected-map residual, so the
factor cancels from relative decay tests. At a completed iteration $t$,
restart when any of the following holds:
\begin{equation}\label{eq:restart-rules}
\begin{aligned}
 &\varepsilon_{r,t}\leq\alpha_1\varepsilon_{r,1}
       &&\text{(sufficient decay)},\\
 &\varepsilon_{r,t}\leq\alpha_2\varepsilon_{r,1},
   \quad\varepsilon_{r,t}>\varepsilon_{r,t-1}
       &&\text{(decay followed by an increase)},\\
 &t\geq\alpha_3 k
       &&\text{(long cycle)}.
\end{aligned}
\end{equation}
Here $0<\alpha_1<\alpha_2<1$ and $0<\alpha_3<1$; the first two tests
are applied for $t\geq2$. This indexing uses only completed sweeps.
The new anchor is the last sweep output $\bar w^{r,\tau_r}$.
Related restart strategies for LP use normalized duality gaps or weighted
KKT residuals~\cite{applegate2021practical,lu2023cupdlp,lu2023cupdlp_c}.

\subsection{Penalty updates}\label{subsec:sigma}
The metric in~\eqref{def:M} satisfies
\begin{equation}\label{eq:metric-expansion}
\begin{aligned}
 \|w-w^*\|_{\mathcal M_\sigma}^2
 &=\sigma\bigl(\|y-y^*\|_{\mathcal T_1}^2+\|A^*(y-y^*)\|^2\bigr)\\
 &\quad+2\langle A^*(y-y^*),x-x^*\rangle
  +\frac{\|x-x^*\|^2}{\sigma}.
\end{aligned}
\end{equation}
When both coefficients of $\sigma$ and $\sigma^{-1}$ are positive,
minimizing this expression over $\sigma>0$ gives
\begin{equation}\label{eq:sigma}
 \sigma_{\rm ideal}
 =\frac{\|x-x^*\|}
 {\sqrt{\|y-y^*\|_{\mathcal T_1}^2+\|A^*(y-y^*)\|^2}}.
\end{equation}
This minimizes the anchor-distance factor, not the entire KKT bound,
whose prefactor also depends on $\sigma$.

Since the solution is unknown, HPR-LP replaces solution differences by
progress over the completed cycle~\cite[Section~3.2]{chen2024hpr}:
\begin{equation}\label{eq:Dx&Dy}
\begin{aligned}
 \Delta_x&=\|\bar x^{r,\tau_r}-x^{r,0}\|,\\
 \Delta_y&=\bigl(\|\bar y^{r,\tau_r}-y^{r,0}\|_{\mathcal T_1}^2
       +\|A^*(\bar y^{r,\tau_r}-y^{r,0})\|^2\bigr)^{1/2}.
\end{aligned}
\end{equation}
The candidate update is $\sigma_{\rm cand}=\Delta_x/\Delta_y$.
For the two proximal choices in Section~\ref{subsec:projection}, this becomes
\begin{equation}\label{eq:sigma-cases}
 \sigma_{\rm cand}=
 \begin{cases}
 \displaystyle\frac{\|\bar x^{r,\tau_r}-x^{r,0}\|}
 {\|A^*(\bar y^{r,\tau_r}-y^{r,0})\|},&\mathcal T_1=0,\\[8pt]
 \displaystyle\frac{1}{\sqrt{\lambda_A}}
 \frac{\|\bar x^{r,\tau_r}-x^{r,0}\|}
 {\|\bar y^{r,\tau_r}-y^{r,0}\|},
 &\mathcal T_1=\lambda_A I_m-AA^*.
 \end{cases}
\end{equation}
The second rule is closely related to the primal-weight update in
PDLP~\cite[Algorithm~3]{applegate2021practical}.
The candidate requires safeguards against zero or unreliable displacements
and extreme changes, as described in~\cite[Section~3.2]{chen2024hpr};
related smoothing strategies appear in~\cite{chen2025hpr}.
An undefined or inadmissible candidate must not replace the current
positive penalty.

\subsection{Termination}
Restart residuals measure progress in the current metric. Termination
instead requires the prescribed primal feasibility, dual feasibility, and
primal--dual gap checks; the original HPR-LP criteria are given
in~\cite[Section~4.1]{chen2024hpr}. When preprocessing or scaling is used,
the recovered candidate must be assessed in the original model.

\section{GPU preprocessing}\label{sec:preprocessing}

\href{https://github.com/PolyU-IOR/GPU-Preprocessor}{GPU preprocessing} reduces the LP before HPR iteration through folding
and presolve. The underlying framework, correctness results, and
experiments are presented in~\cite{chen2026gpu_preprocessing}. We
summarize the main components and describe implementation extensions
for HPR-LP.

\subsection{Workflow and scheduling}
The workflow starts with structural reductions, then applies folding,
followed by four presolve phases. Startup precedes folding because
aggregation can obscure patterns needed by the structural rules.
Cleanup removes basic redundancies; Propagation tightens bounds and
performs dual fixing; Aggregation eliminates variables by substitution;
Parallel Detection finds proportional rows and columns.
Algorithm~\ref{alg:preprocessing-workflow} prioritizes the first two,
less expensive phases and returns to Cleanup when a later phase
enables further reductions.

\begin{algorithm}[htbp]
\caption{GPU preprocessing schedule, adapted from~\cite{chen2026gpu_preprocessing}}
\label{alg:preprocessing-workflow}
\small
\begin{algorithmic}[1]
\REQUIRE LP $P$ and positive phase budget $N_{\rm phase}$.
\STATE $P\gets\textsc{Startup}(P)$; set $q=1$.
\STATE $P\gets\textsc{Folding}(P)$.
\WHILE{true}
 \STATE \textbf{Cleanup:} apply \textsc{Cleanup} until $P$ is unchanged.
 \STATE \textbf{break if} $q\geq N_{\rm phase}$.
 \STATE $P\gets\textsc{Propagation}(P)$; $q\gets q+1$.
 \IF{Propagation sufficiently reduces $\operatorname{nnz}(A)$ or tightens bounds}
  \STATE \textbf{go to Cleanup}.
 \ENDIF
 \STATE \textbf{break if} $q\geq N_{\rm phase}$.
 \STATE $P\gets\textsc{Aggregation}(P)$; $q\gets q+1$.
 \IF{Aggregation sufficiently reduces $\operatorname{nnz}(A)$}
  \STATE \textbf{go to Cleanup}.
 \ENDIF
 \STATE \textbf{break if} $q\geq N_{\rm phase}$.
 \STATE $P\gets\textsc{Parallel Detection}(P)$; $q\gets q+1$.
 \IF{Parallel Detection sufficiently reduces $\operatorname{nnz}(A)$}
  \STATE \textbf{go to Cleanup}.
 \ENDIF
 \STATE \textbf{break}.
\ENDWHILE
\ENSURE The preprocessed LP $P$.
\end{algorithmic}
\end{algorithm}

Here $\operatorname{nnz}(A)$ is the number of matrix nonzeros. The
budget counts Startup and subsequent non-Cleanup phase executions;
folding and Cleanup are excluded. The returned model is passed to
scaling and HPR-LP. Recovery is a separate stage.

\subsection{Folding}
Folding aggregates rows and columns using a stable partition pair.
Within each row block, row bounds agree and coefficient sums into each
column block agree. The analogous conditions hold for columns, including
their costs and variable bounds. These conditions permit a smaller
folded LP with a solution-recovery map; see~\cite[Section~2.1]{chen2026gpu_preprocessing}.

The GPU implementation starts by grouping rows with equal
$((l_c)_i,(u_c)_i)$ and columns with equal
$(c_j,(l_v)_j,(u_v)_j)$. It then alternates row and column refinement.
Independent Gaussian weights are assigned to the current blocks and
expanded to vectors. Products with $A$ and $A^*$ produce signatures
encoding coefficient sums into the opposite partition.
Parallel sorting groups equal signatures within existing blocks, and
prefix scans assign new block indices.
Refinement stops when the partitions stabilize, or returns singleton
partitions when both partition sizes exceed a prescribed fraction of
the original dimensions.
The probability-one guarantee and refinement-round bound are given
in~\cite[Theorem~3]{chen2026gpu_preprocessing}.

\subsection{Presolve and implementation extensions}
Our implementation extends Startup with equality-row
sparsification, Poisson elimination, and two-row variable elimination.
Equality-row sparsification cancels common coefficients through
certified row differences without changing the variables or their
bounds. For certified three-dimensional Dirichlet Poisson stencils,
cuFFT-based Green responses eliminate grid-state variables and form
observation equations in the source variables. Other structural rules
include column domination, zero-cost product variables, and
piecewise-linear reformulation.

\begin{figure}[htbp]
\centering
\includegraphics[width=0.88\linewidth]{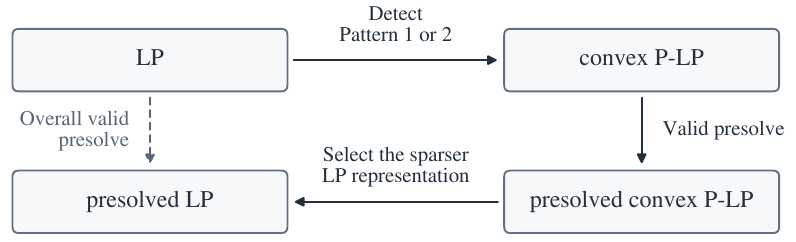}
\caption{Piecewise-linear reformulation workflow, reproduced
from~\cite[Figure~1]{chen2026gpu_preprocessing}.}
\label{fig:pwl-reformulation}
\end{figure}

The piecewise-linear rule identifies an underlying convex
piecewise-linear program, simplifies it, and chooses a sparser LP
representation (Figure~\ref{fig:pwl-reformulation}).
Its validity and proof are given in~\cite[Theorem~4 and Appendix~B]{chen2026gpu_preprocessing}.

The Scan--Plan--Apply design shares sparse-data scans across rules.
Scan collects row and column information while respecting deletion masks;
Plan accumulates compatible deletions and vector updates; Apply commits
the updates and rebuilds $A$ and $A^*$ once per batch.
For doubleton equations, a subset with disjoint variable sets avoids
conflicting bound and objective updates. Atomic operations accumulate
contributions when different substitutions affect the same row.
For proportional rows and columns, randomized SpMV screening and
hash-based grouping restrict explicit comparisons to candidate groups.
These implementation mechanisms are described
in~\cite[Section~3.3]{chen2026gpu_preprocessing}.

\section{Diagonal scaling}\label{sec:scaling}

Scaling changes the coordinates in which HPR-LP operates. We describe
a Curtis--Reid (CR) pipeline and extended root-mean-square (RMS)
scaling.

\subsection{Coordinate transformations}\label{scaling:coordinates}
For positive diagonal matrices $D_r$ and $D_c$, transform
the data in~\eqref{model:primalLP} as
\begin{equation}\label{scaling:data-transform}
\begin{aligned}
 \widetilde A&=D_rAD_c,&\widetilde c&=D_cc,\\
 (\widetilde l_c,\widetilde u_c)&=(D_rl_c,D_ru_c),&
 (\widetilde l_v,\widetilde u_v)&=(D_c^{-1}l_v,D_c^{-1}u_v).
\end{aligned}
\end{equation}
Recover the variables before this diagonal transformation by
\begin{equation}\label{scaling:recovery}
 x=D_c\widetilde x,\qquad
 y=D_r\widetilde y,\qquad
 z=D_c^{-1}\widetilde z.
\end{equation}
These invertible transformations preserve feasibility and optimality.
In particular,
$\widetilde A^*\widetilde y+\widetilde z=\widetilde c$
is equivalent to $A^*y+z=c$.

The standard pipeline considered here is
\begin{equation}\label{scaling:baseline}
 \mathrm{CR}(20)\ \longrightarrow\ \mathrm{Ruiz}(10)\
 \longrightarrow\ \mathrm{PC}\ \longrightarrow\ bc.
\end{equation}
The integers denote complete sweeps. Ruiz scaling balances row and
column infinity norms~\cite{scaling-ruiz2001}; PC denotes
Pock--Chambolle scaling~\cite{scaling-pock2011}. The final $bc$ step
normalizes bounds and objective coefficients. Every matrix scaling
also transforms the LP vectors consistently.

\subsection{Curtis--Reid scaling on the GPU}\label{scaling:cr}
CR scaling minimizes the squared logarithms of nonzero matrix
magnitudes~\cite{scaling-curtis1972}. Set
$\Omega=\{(i,j):A_{ij}\ne0\}$ and $\ell_{ij}=\log|A_{ij}|$.
With row and column logarithmic multipliers $a\in\mathbb R^m$
and $b\in\mathbb R^n$, the model is
\begin{equation}\label{scaling:cr-objective}
 \min_{a,b}\ \frac12\sum_{(i,j)\in\Omega}(\ell_{ij}+a_i+b_j)^2.
\end{equation}
The original CR algorithm uses conjugate gradients. The HPR-LP
implementation instead performs alternating least-squares (ALS)
updates. If $\delta_i$ and $\nu_j$ are the numbers of nonzeros in row
$i$ and column $j$, respectively, one sweep is
\begin{equation}\label{scaling:als}
\begin{aligned}
 a_i^{s+1}&=-\frac{1}{\delta_i}
       \sum_{j:(i,j)\in\Omega}(\ell_{ij}+b_j^s),\\
 b_j^{s+1}&=-\frac{1}{\nu_j}
       \sum_{i:(i,j)\in\Omega}(\ell_{ij}+a_i^{s+1}).
\end{aligned}
\end{equation}
Empty rows and columns retain the identity factor.
Starting from $a^0=b^0=0$, 20 sweeps give
$D_r=\operatorname{diag}(\exp(a^{20}))$ and
$D_c=\operatorname{diag}(\exp(b^{20}))$.
Each block update exactly minimizes~\eqref{scaling:cr-objective} in that
block before numerical safeguards are applied.

ALS is block Gauss--Seidel on the CR normal equations. Adding a constant
to all row multipliers and subtracting it from all column multipliers
in a connected support component leaves the scaled matrix unchanged.
This gauge freedom explains why different factor pairs may represent
the same scaled matrix. ALS and conjugate gradients target the same
least-squares objective, but finite sweeps, gauges, and safeguards can
give different factors.

One ALS sweep visits the entries of $A$ and $A^*$ and costs
$O(\operatorname{nnz}(A))$. The implementation caches
$-\log(\max\{|A_{ij}|,10^{-300}\})$, avoiding repeated logarithm
evaluation. Rows are grouped by degree to reduce GPU workload imbalance;
long rows can use cooperative loading with ordered accumulation.
The row and column phases remain sequential, although updates within
each phase are parallel. Factor-spread checks and bounded exponentiation
can reject an unsafe candidate before the LP is changed.

\subsection{Extended RMS scaling}\label{scaling:rms}
CR scaling uses matrix entries alone. Extended RMS also incorporates
constraint bounds and objective coefficients when choosing factors.
For a vector with $N>0$ nonzeros,
$\operatorname{RMS}(v)=\|v\|/\sqrt N$.
Define a row-bound magnitude vector by
\begin{equation}\label{scaling:bound-magnitude}
 \widehat b_i=\max\{\phi((l_c)_i),\phi((u_c)_i)\},\qquad
 \phi(q)=
 \begin{cases}|q|,&q\text{ finite},\\0,&\text{otherwise}.\end{cases}
\end{equation}
The factors are computed from the virtual symmetric matrix
\begin{equation}\label{scaling:augmented-matrix}
 H=\begin{bmatrix}
 0&A^*&c\\ A&0&\widehat b\\c^*&\widehat b^*&0
 \end{bmatrix}.
\end{equation}
Its entries are accessed through the existing GPU arrays; $H$ is not
assembled and introduces no additional optimization variables.

Starting from $q^0=\mathbf1$, write
$H^s=\operatorname{diag}(q^s)^{-1}H\operatorname{diag}(q^s)^{-1}$.
For a nonempty row with support $\mathcal N_i=\{j:H_{ij}\ne0\}$, set
\begin{equation}\label{scaling:augmented-update}
 \rho_i^s=\left(\frac{1}{|\mathcal N_i|}
                  \sum_{j\in\mathcal N_i}|H_{ij}^s|^2\right)^{1/2},
 \qquad q_i^{s+1}=q_i^s\sqrt{\rho_i^s},
 \qquad s=0,\ldots,9.
\end{equation}
All updates in a round use $q^s$. The implementation uses logarithmic
coordinates, bounds accumulated divisors within $[10^{-8},10^8]$,
and uses identity updates for empty rows and numerically negligible
cases. After ten rounds, apply~\eqref{scaling:data-transform} with
\begin{equation}\label{scaling:rms-factors}
 D_c=\operatorname{diag}(q_1^{10},\ldots,q_n^{10})^{-1},\qquad
 D_r=\operatorname{diag}(q_{n+1}^{10},\ldots,q_{n+m}^{10})^{-1}.
\end{equation}
The auxiliary divisor $q_{n+m+1}^{10}$ influences factor selection
but is not directly applied to the LP. The resulting positive
diagonal transformation preserves the LP regardless of whether the
ten-round balancing procedure has converged.

\subsection{Final normalization}\label{scaling:pc-bc}
For the current scaled matrix $\widetilde A$, PC uses row and column
divisors computed from the same matrix:
\begin{equation}\label{scaling:pc}
 p_i=\left(\sum_j|\widetilde A_{ij}|\right)^{1/2},\qquad
 d_j=\left(\sum_i|\widetilde A_{ij}|\right)^{1/2}.
\end{equation}
Apply~\eqref{scaling:data-transform} with
$D_r=\operatorname{diag}(p)^{-1}$ and
$D_c=\operatorname{diag}(d)^{-1}$, replacing zero or extremely small
divisors by one.

Using the data after PC, let
$\beta_b=1+\|\widehat b_{\rm P}\|$ and
$\beta_c=1+\|c_{\rm P}\|$, where $\widehat b_{\rm P}$ is defined
by~\eqref{scaling:bound-magnitude}.
The $bc$ step leaves the matrix unchanged, divides all row and variable
bounds by $\beta_b$, and divides the objective vector by $\beta_c$.
Recovery to the coordinates immediately before this step is
\begin{equation}\label{scaling:bc-recovery}
 x_{\rm P}=\beta_b x_{\rm final},\qquad
 y_{\rm P}=\beta_c y_{\rm final},\qquad
 z_{\rm P}=\beta_c z_{\rm final}.
\end{equation}
Objective values acquire the factor $\beta_b\beta_c$.
Apply this recovery first, then undo each diagonal transformation
in reverse order using~\eqref{scaling:recovery}. Any preprocessing
recovery follows before original-model checks.

\section{ISpMV: reusable execution plans}\label{sec:spmv}

Each HPR sweep~\eqref{HPR-LP-2} applies the same operators to changing
vectors. Matrix- and hardware-dependent preparation can therefore be
amortized over the iterations~\cite{spmv-oski2005}.
ISpMV~\cite{spmv-ispmv} constructs a reusable plan
$\mathcal P_{B,\mathcal H}=(R,E,S)$ for each $B\in\{A,A^*\}$ and GPU
resource profile $\mathcal H$. Its components specify an algebraic
representation $R$, storage encoding $E$, and parallel schedule $S$.
For conforming vectors $u,v$ and scalars $\alpha,\beta\in\mathbb R$,
its exact-arithmetic action is
\begin{equation}\label{spmv:apply}
 v_{\rm new}=\alpha Bu+\beta v
 =\operatorname{Apply}(\mathcal P_{B,\mathcal H};u,v,\alpha,\beta).
\end{equation}

\paragraph{Algebraic representation}
Our GPU detectors verify supported row, block, and incidence patterns
against the complete support and coefficients of $B$, and construct both
an exact decomposition and its input/output maps:
\begin{equation}\label{spmv:representation}
 B=B_{\rm res}+\sum_{\ell=1}^{L}U_\ell C_\ell V_\ell,
 \qquad
 Bu=B_{\rm res}u+\sum_{\ell=1}^{L}U_\ell\bigl[C_\ell(V_\ell u)\bigr].
\end{equation}
Here $V_\ell$ gathers or aggregates inputs, $C_\ell$ applies a core
operator, and $U_\ell$ maps its outputs to matrix rows. The factors are
implemented through stored cores, index maps, and implicit formulas.
Proportional rows reuse representative dot products, and shared weighted
group sums are evaluated once; both exploit common-subexpression
reuse~\cite{spmv-oppelstrup2013}. Repeated blocks share coefficient cores,
while regular incidence patterns admit implicit index formulas.
Uncovered entries remain in $B_{\rm res}$; $L=0$ recovers conventional
sparse execution.

\paragraph{Storage encoding}
Building on value and index compression~\cite{spmv-kourtis2010}, we
construct coefficient dictionaries, local column offsets, and reduced-width
coefficient storage on the GPU. Encoding eligibility is checked against the
values and index ranges required by $R$, subject to exact recovery of the
represented values and indices. Each admitted encoding is paired with a
compatible kernel that consumes it directly with FP64 arithmetic, without
materializing a decoded matrix.

\paragraph{Parallel scheduling}
We use thread- and warp-level row processing~\cite{spmv-bell2009} and
merge-based row/nonzero partitioning~\cite{spmv-merrill2016} as scheduling
primitives. ISpMV constructs row groups, long-row segments, and reduction
descriptors on the GPU for compatible representations and encodings.
Structured components and the sparse residual can use different schedules
within one plan, subject to the dependencies in~\eqref{spmv:representation}
and the GPU resource limits.

Numerical probes screen the FP64 evaluation separately from algebraic
verification. Bounded trials compare the cost of complete plans, including
mappings, residual products, and reductions. The selected plan is reused
while $B$ remains unchanged.

\section{A sieving technique}\label{sec:sieving}

Sieving adaptively fixes selected variables at bounds and temporarily
removes selected constraints to reduce the cost of SpMV operations per
iteration. Within each restart cycle of Algorithm~\ref{alg:HPR_LP},
we distinguish full and sieved iterations.
Full iterations use the full working LP, after any preprocessing and
scaling; sieved iterations incorporate the selected variable fixings
and constraint removals.
Unlike presolve, these selections are reconsidered at full iterations
and reset at each restart.

\subsection{Full and sieved iterations}
At the start of cycle $r$, initialize the constraint and variable
sieving index sets as
\begin{equation}\label{eq:sieving-initial}
 \mathcal I_0=\{i:y_i^{r,0}=0\},\;
 \mathcal J_0=\{j:x_j^{r,0}\in\{(l_v)_j,(u_v)_j\}\}.
\end{equation}
Full iterations occur at indices $0=t_0<t_1<\cdots$.
At $t=t_h$, compute
$\bar w^{r,t_h+1}=\mathcal D_{\sigma_r}(w^{r,t_h})$ and retain only
candidates whose sweep outputs agree with their anchor values:
\begin{equation}\label{eq:sieving-selection}
 \mathcal I_{h+1}
 =\{i\in\mathcal I_h:\bar y_i^{r,t_h+1}=y_i^{r,0}\},\;
 \mathcal J_{h+1}
 =\{j\in\mathcal J_h:\bar x_j^{r,t_h+1}=x_j^{r,0}\}.
\end{equation}
We then define
\begin{equation}\label{eq:sieving-sets}
\begin{aligned}
 \mathcal K_{h+1}
 &=\{s\in\mathbb R^m:(l_c)_i\leq s_i\leq(u_c)_i
                      \text{ for }i\notin\mathcal I_{h+1}\},\\
 \mathcal C_{h+1}
 &=\{x\in\mathbb R^n:(l_v)_j\leq x_j\leq(u_v)_j
       \text{ for }j\notin\mathcal J_{h+1},\;
       x_j=x_j^{r,0}\text{ for }j\in\mathcal J_{h+1}\}.
\end{aligned}
\end{equation}
For $h\geq1$, let $\mathcal D_{\sigma_r}^{h}$
denote the sweep in Algorithm~\ref{alg:sp-HPR} with $\mathcal K$ and
$\mathcal C$ replaced by $\mathcal K_h$ and $\mathcal C_h$ in the
augmented Lagrangian~\eqref{def:augmented-lagrangian}.
For $t_{h-1}<t<t_h$, a sieved iteration computes
$\bar w^{r,t+1}=\mathcal D_{\sigma_r}^{h}(w^{r,t})$.
In both cases, we then apply the reflection and Halpern steps in
Algorithm~\ref{alg:HPR_LP} to obtain $w^{r,t+1}$.

In the projection-based implementation, sieving reduces SpMV work in
sieved iterations by omitting rows corresponding to removed constraints
and reusing cached contributions from fixed variables.

\subsection{Scheduling and restart}
Let $0\leq\theta<1$ and let $p\geq1$ be an integer.
At a full iteration $t_h$, test
\begin{equation}\label{eq:sieving-controller}
 \bigl\langle w^{r,t_h}-\bar w^{r,t_h+1},
       \mathcal M_r(w^{r,t_h}-w^{r,0})\bigr\rangle
 \leq-(1-\theta)t_h
       \|w^{r,t_h}-\bar w^{r,t_h+1}\|_{\mathcal M_r}^2.
\end{equation}
If the test holds, the next full iteration is $t_{h+1}=t_h+p$;
otherwise it is $t_{h+1}=t_h+1$. Selection~\eqref{eq:sieving-selection}
is performed at every full iteration, independently of this test.
Algorithm~\ref{alg:sieving} replaces one restart cycle of
Algorithm~\ref{alg:HPR_LP}.

\begin{algorithm}[htbp]
\caption{One HPR-LP cycle with sieving}
\label{alg:sieving}
\small
\begin{algorithmic}[1]
\REQUIRE Anchor $w^{r,0}$, penalty $\sigma_r$, $\theta\in[0,1)$,
and integer $p\geq1$.
\STATE Set $t=0$, $h=0$, $t_0=0$, and initialize the sets
in~\eqref{eq:sieving-initial}.
\LOOP
 \STATE Set $\mathrm{full}=(t=t_h)$.
 \IF{$\mathrm{full}$}
  \STATE $\bar w^{r,t+1}=\mathcal D_{\sigma_r}(w^{r,t})$.
  \STATE Set $t_{h+1}=t_h+p$ if~\eqref{eq:sieving-controller} holds,
         and $t_{h+1}=t_h+1$ otherwise.
  \STATE Update the sets by~\eqref{eq:sieving-selection}--\eqref{eq:sieving-sets};
         set $h=h+1$.
 \ELSE
  \STATE $\bar w^{r,t+1}
       =\mathcal D_{\sigma_r}^{h}(w^{r,t})$.
 \ENDIF
 \STATE $\hat w^{r,t+1}=2\bar w^{r,t+1}-w^{r,t}$.
 \STATE $\displaystyle w^{r,t+1}
       =\frac{1}{t+2}w^{r,0}+\frac{t+1}{t+2}\hat w^{r,t+1}$.
 \STATE Set $t=t+1$ and increment the total iteration counter.
 \IF{$\mathrm{full}$ and the original-model termination test holds}
  \STATE \textbf{return} $\bar w^{r,t}$ as the final candidate.
 \ENDIF
 \IF{$\mathrm{full}$ and a restart criterion holds}
  \STATE Set $\tau_r=t$ and update $\sigma_{r+1}$ as in Section~\ref{subsec:sigma}.
  \STATE \textbf{return} $\bar w^{r,t}$ as the next anchor.
 \ENDIF
\ENDLOOP
\end{algorithmic}
\end{algorithm}

Restart residuals are evaluated at full sweeps in the full working
metric $\mathcal M_r$.
The decay reference is the first full sweep of the cycle; the
local-increase test compares successive full-sweep residuals.
Only full iterations trigger restart or termination.
With $p=1$, every iteration is full and the scheme reduces to the
restarted base method.

\FloatBarrier
\section{Multi-GPU parallelization}\label{sec:multigpu}

\begin{figure}[!htbp]
 \centering
 \includegraphics[width=\linewidth]{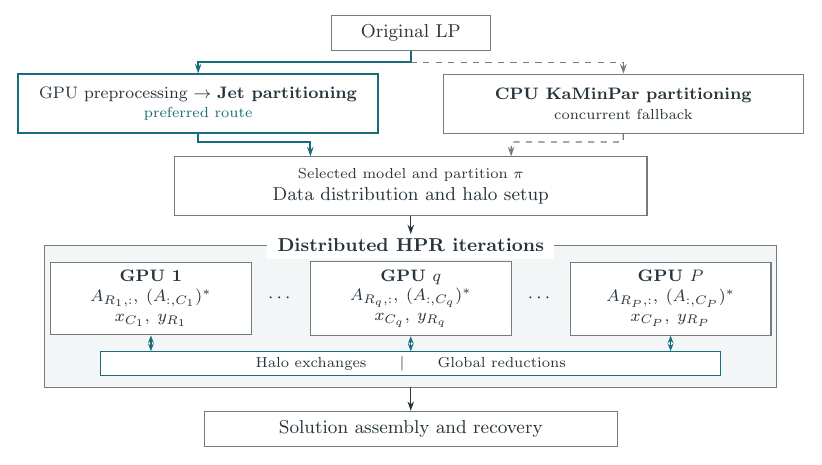}
 \caption{Multi-GPU parallelization of HPR-LP. The selected
 preparation path determines the working model and its row--column
 distribution; local HPR updates communicate through halo exchanges
 and global reductions.}
 \label{fig:multigpu-architecture}
\end{figure}

We extend HPR-LP to multiple GPUs to accelerate large-scale LPs
and to handle instances that exceed the memory capacity of a single
GPU. A joint row--column partition assigns ownership of the primal
and dual variables to $P$ GPUs, determining the local matrix data
and the communication required in each HPR sweep. We choose the
partition to balance computational work while limiting communication.
Figure~\ref{fig:multigpu-architecture} outlines the overall architecture.

\FloatBarrier

\subsection{Joint partitioning model}
Let $G=(V,E)$ be the bipartite graph of the working matrix
$A\in\mathbb R^{m\times n}$, with row vertices $v_i^r$, column vertices
$v_j^c$, and an edge $(v_i^r,v_j^c)$ for each stored entry $A_{ij}$.
A partition $\pi:V\to\{1,\ldots,P\}$ defines
$R_q=\{i:\pi(v_i^r)=q\}$ and $C_q=\{j:\pi(v_j^c)=q\}$.
As in cuOpt mPDLP~\cite{multigpu-mpdlp2026}, GPU $q$ owns
$y_{R_q}$ and $x_{C_q}$ and computes $(Ax)_{R_q}$ and $(A^*y)_{C_q}$.
To balance the work of these products while limiting the exchange
of remote inputs, we use the partitioning model
\begin{equation}\label{eq:multigpu-partition}
 \begin{aligned}
 \min_{\pi}\quad&\operatorname{cut}(\pi)
 :=\sum_{(v_i^r,v_j^c)\in E}\mathbf 1_{\{\pi(v_i^r)\ne\pi(v_j^c)\}}\\
 \text{subject to}\quad& W_q(\pi):=\sum_{v:\pi(v)=q}\omega_v\leq B,
 \qquad q=1,\ldots,P.
 \end{aligned}
\end{equation}
Here $\omega_v\geq0$ are prescribed vertex weights and $B$ is a
per-GPU workload budget. The constraint approximates computational
balance; the objective minimizes an upper bound on halo communication volume.

\begin{samepage}
The unit weights $\omega_v=1$ used in cuOpt's
partitioner~\cite{multigpu-cuopt-source2026} balance the number of
owned coordinates, $W_q(\pi)=|R_q|+|C_q|$. We use
$\omega_v=\deg(v)$, giving
\begin{equation}\label{eq:multigpu-load}
 W_q(\pi)=\operatorname{nnz}(A_{R_q,:})+
          \operatorname{nnz}(A_{:,C_q}).
\end{equation}
This choice accounts for heterogeneous row and column densities.
The constraint approximates computational load balance, since it
aggregates two sequential products whose per-GPU workloads and
kernel efficiencies may differ.
\end{samepage}

\begin{samepage}
To relate the objective to communication, consider the remote inputs
required under the same ownership. Counting each coordinate once
per destination gives the input halos and aggregate payload per
product pair
\begin{equation}\label{eq:multigpu-halos}
 \begin{aligned}
 H_q^x&:=\{j\notin C_q:\ \exists\,i\in R_q\text{ with }(v_i^r,v_j^c)\in E\},\\
 H_q^y&:=\{i\notin R_q:\ \exists\,j\in C_q\text{ with }(v_i^r,v_j^c)\in E\},\\
 \mathcal V(\pi)&:=\sum_{q=1}^P\bigl(|H_q^x|+|H_q^y|\bigr)
                  \leq2\operatorname{cut}(\pi).
 \end{aligned}
\end{equation}
\end{samepage}
With uniform cost per transmitted scalar, $\mathcal V(\pi)$
approximates aggregate halo communication cost. The inequality
follows because each cut edge induces at most one required scalar
input per product. For the LP partitions considered here, the bound
is typically strict because several cut edges share a remote
coordinate that is transmitted only once to the receiving GPU.

\subsection{Distributed implementation}
By default, GPU preprocessing is followed by Jet
partitioning~\cite{multigpu-jet2024} of the reduced matrix, while CPU
\mbox{KaMinPar}~\cite{multigpu-kaminpar2021} concurrently partitions
the original matrix as a fallback. Both partitioners
solve~\eqref{eq:multigpu-partition} heuristically on their respective
working matrices. The selected working model and partition determine
the distributed operators and halo exchanges used for spectral
estimation and subsequent HPR iterations.

A full sweep of~\eqref{HPR-LP-2} first exchanges $y_{H_q^y}$ for the
primal update and then $(2\bar x-x)_{H_q^x}$ for the dual update.
Unique output ownership avoids inter-GPU reductions of partial
outputs. Projections, reflections, and Halpern updates are local;
global reductions supply convergence measures and the quantities
used for restart and penalty updates (Section~\ref{sec:implementation}).
Sieving (Section~\ref{sec:sieving}) updates local operators and halo
dependencies while preserving ownership. At termination, the
distributed solution is assembled and mapped back to the original model.
\FloatBarrier

\section{OptBin: model storage and loading}\label{sec:optbin}

Large LPs incur costs before iteration begins: reading source files,
decoding values, and constructing sparse arrays.
OptBin addresses this stage through a numerical schema, lossless block
encodings, and a native CPU reader~{\cite{optbin-software}}.
This section summarizes its LP-relevant design and reported loading measurements.

\subsection{Numerical model and storage}
OptBin stores the matrix, objective coefficients, and row and variable
bounds of~\eqref{model:primalLP}, together with the objective sense and
any constant.
Losslessness means exact reconstruction of the canonical binary64
numerical arrays after declared parsing and canonicalization; it does
not mean reproduction of decimal tokens, names, or text ordering.
Duplicate summation and coordinate sorting, when used, belong to
canonicalization. Bounds may be infinite; NaNs and nonfinite ordinary
coefficients are rejected.

In compressed sparse column (CSC) form, a matrix with $n$ columns
and $N_A$ stored entries occupies
\begin{equation}\label{optbin:footprint}
 8(n+1)+w_{\rm ind}N_A+8N_A\quad\text{bytes}
\end{equation}
in the decoded layout, excluding vectors and metadata.
Pointers are 64-bit offsets, whereas row indices use
$w_{\rm ind}=4$ or $8$ bytes according to the row-dimension limit.
Entry counts and row dimensions impose different width requirements.
For example, the stored header of Hinder's
\texttt{prod\_200\_500\_02} has $4{,}443{,}349{,}500$ entries
and $25{,}476{,}450$ rows~\cite{optbin-hinder,optbin-hinder-generators}.
It needs 64-bit pointers but admits 32-bit row indices.
Widening those indices would add 16.55 GiB by~\eqref{optbin:footprint}.

Integer blocks use raw values or difference encodings.
Coefficient blocks choose among raw binary64, repeated bit patterns,
exact narrow integers, exact binary32 values, and dictionary indices.
Narrowing is allowed only when decoding reproduces the binary64 bits.
Optional Zstandard compression~\cite{optbin-zstd} follows block encoding.
The writer selects the smallest encoding before compression,
so its default does not guarantee the smallest compressed file.

The container records types, shapes, offsets, codec identifiers, and
checksums over decoded data. The reader validates sparse ordering,
coverage, and cross-array consistency. Both OptBin and the encoded
HDF5 baseline~\cite{optbin-hdf5} use the same block decoder and validator.
Full loading still allocates the expanded numerical arrays: smaller
files do not imply smaller solver working memory.

\subsection{Large-LP loading measurements}
The two LPs below come from Hinder's benchmark collection and generators
associated with the PDLP study~\cite{optbin-hinder,optbin-hinder-generators,optbin-pdlp-large}.
We abbreviate \texttt{prod\_200\_400\_02} as \texttt{prod} and
\texttt{mcf\_5000\_100\_1000} as \texttt{mcf}.
Their dimensions are
\[
\begin{array}{lrrr}
 &m&n&\operatorname{nnz}(A)\\
 \texttt{prod}&16{,}361{,}200&64{,}560{,}800&2{,}310{,}559{,}600\\
 \texttt{mcf}&5{,}525{,}100&502{,}500{,}100&1{,}007{,}500{,}100.
\end{array}
\]
\paragraph{Practical compressed-file paths}
Table~\ref{optbin:tab:public-loading} compares original MPS.gz loading
through Gurobi 13.0.3 with the native binary paths.
Both LPs were validated through independent source-model checks and
decoded-array comparisons. The binary paths return validated numerical arrays;
Gurobi returns its solver model, so the endpoints have different
internal representations.
Loading times were measured on an Intel Xeon Gold 6548Y+ shared host,
using warm files and three shuffled fresh-process trials per format.
Measurements for the two LPs were collected concurrently, using serial
native decoding and one Gurobi thread.
Environment startup and post-load comparison were excluded.

\begin{table}[htbp]
\centering\small
\begin{tabular}{lrrrr}
\toprule
&\multicolumn{2}{c}{\texttt{prod}}&\multicolumn{2}{c}{\texttt{mcf}}\\
Representation&GiB&Seconds&GiB&Seconds\\
\midrule
MPS.gz (Gurobi)&7.294&342.27&10.553&386.89\\
Raw HDF5.gz&0.533&101.39&4.301&117.72\\
OptBin&0.441&53.13&4.460&53.59\\
OptBin.gz&0.421&55.41&4.314&73.37\\
\bottomrule
\end{tabular}
\caption{Stored size and median loading time.
GiB denotes $2^{30}$ bytes. Binary gzip paths include decompression to
a seekable temporary file and native loading.}
\label{optbin:tab:public-loading}
\end{table}

These complete loader paths favor OptBin in time on the two cases,
but raw HDF5.gz is smaller on \texttt{mcf}.
OptBin peak process RSS was 28.00 and 26.32 GiB, respectively;
the other binary paths were within 0.01 GiB.
Gurobi read-phase RSS was 57.67 and 70.41 GiB and included its model
and metadata.
Conversion, gzip preparation, solver setup, preprocessing, GPU upload,
and optimization were outside the reported loading times.

\paragraph{Matched-codec container comparison}
Table~\ref{optbin:tab:matched-loading} compares OptBin and HDF5 with
identical encoded block bytes, alongside typed raw HDF5 and a chunked
shuffle/gzip baseline. Both encoded paths use the same numerical
decoder, allocations, checksum checks, and full validator.

\begin{table}[htbp]
\centering\small
\begin{tabular}{lrrrr}
\toprule
&\multicolumn{2}{c}{\texttt{prod}}&\multicolumn{2}{c}{\texttt{mcf}}\\
Representation&GiB&Seconds&GiB&Seconds\\
\midrule
OptBin&0.441&51.19&4.460&52.16\\
Encoded HDF5&0.441&51.13&4.460&52.32\\
Raw HDF5&27.992&42.26&26.319&39.00\\
HDF5 shuffle/gzip&0.783&74.51&4.732&101.21\\
\bottomrule
\end{tabular}
\caption{Stored size and median loading time for the container comparison.
Times are medians of three shuffled fresh-process trials on warm files,
including reading, decoding, checksums, and full validation.
Fingerprinting and explicit byte comparisons are outside timing.}
\label{optbin:tab:matched-loading}
\end{table}

The encoded containers have nearly identical sizes and loading medians
on these two LPs. Raw HDF5 loads faster from warm files but uses more
stored bytes. All paths allocate essentially the same decoded arrays;
peak RSS ranges from 27.99 to 28.15 GiB for \texttt{prod} and
26.32 to 26.40 GiB for \texttt{mcf}.
These loading times were measured on the same shared-server processor,
with HDF5 1.14.6, Zstandard 1.4.8, and zlib 1.2.11.

\subsection{Role in an HPR-LP workflow}
OptBin loads the numerical arrays before solver initialization.
The reader is CPU based. Reported process RSS excludes
filesystem page cache, temporary disk files, and later solver allocations.

\begin{samepage}
A binary representation amortizes its conversion cost only through
sufficient reuse. If conversion costs $T_{\rm conv}$ and per-load costs
are $T_{\rm text}>T_{\rm bin}$, then $N_{\rm use}$ uses are beneficial
in conversion-plus-loading time only when
\[
 N_{\rm use}>\frac{T_{\rm conv}}{T_{\rm text}-T_{\rm bin}}.
\]
\end{samepage}

\section{Concluding remarks}\label{sec:conclusion}

HPR-LP combines a semi-proximal HPR iteration with restart and penalty
updates. The projection-based form uses sparse matrix--vector products
and componentwise operations, while the established fixed-parameter
analysis provides nonergodic $O(1/k)$ residual and objective-error bounds.
Its relationships with Halpern PDHG and ergodic PR are stated
at the level of update maps, with explicit parameter and
active-region conditions.

GPU preprocessing and scaling prepare the model. ISpMV provides reusable
execution plans for sparse products. Sieving reduces SpMV work by fixing
variables and removing constraints. Multi-GPU parallelization distributes
sparse products across GPUs, and OptBin addresses model storage and loading.

\bibliographystyle{abbrv}
\bibliography{ref,sections/optbin/references,sections/scaling/references,sections/preprocessing/references,sections/multigpu/references,sections/spmv/references}
\end{document}